\documentclass[sts]{imsart}

\RequirePackage{amsthm,amsmath,amsfonts,amssymb}
\RequirePackage[numbers,sort&compress]{natbib}
\RequirePackage[colorlinks,citecolor=blue,urlcolor=blue]{hyperref}
\RequirePackage{graphicx}
\input preamble.sty
\input RSmacros.sty
\input EPmacros.sty 
\input RStheorems1.sty 
\input RSdebug1.sty 

\usepackage{soul} 
\usepackage{tcolorbox}
\usepackage{adjustbox}
\usepackage{graphicx}
\usepackage{natbib}
\setcitestyle{authoryear,open={(},close={)}}

\startlocaldefs

\def\biglf{\par} 
\theoremstyle{plain}

\theoremstyle{remark}

\endlocaldefs

\begin{document}

\begin{frontmatter}
\title{Paths of Stochastic Processes: \\ a Sudden Turnaround}
\runtitle{Paths: Turnaround}
\begin{aug}
\author[A]{\fnms{Robert}~\snm{Schaback}\ead[label=e1]{schaback@math.uni-goettingen.de}} \and 
\author[B]{\fnms{Emilio}~\snm{Porcu}\ead[label=e2]{emilio.porcu@ku.ac.ae}}


\address[A]{Robert Schaback is Professor Emeritus at Department of Mathematics and Computer Science, Georg-August-Universit\"at  G\"ottingen, Germany\printead[presep={\ }]{e1}.~~~~~~~~~~~~~~~~~~~~~~~~~~~~~~~~~~~~~~~~}
\address[B]{Emilio Porcu is Professor at Department of Mathematics, Khalifa University, and a Senior Fellow at ADIA Lab, Abu Dhabi, The United Arab Emirates\printead[presep={\ }]{e2}.}



\end{aug}

\begin{abstract}
{\em Paths} and their properties ({\em e.g.}, regularity) have a fundamental importance to many theoretical and applied disciplines. Yet, there is a pitfall in the same definition of paths that has been accepted by the literature with no major criticism, nor an attempt to fix it. The commonly accepted definition of paths {starts from a random field} but ignores the problem of setting joint distributions of infinitely many random variables 
{for defining paths properly afterwards}. 
This paper provides a {\em turnaround} that  starts with a given covariance function, then defines paths and {\em finally} a random field.
We show how this approach retains essentially the same properties for Gaussian fields while allowing to construct random fields whose finite dimensional distributions are not Gaussian. Specifically, we start with a kernel $C$ and the associated 
{Reproducing Kernel Hilbert Space} 
${\cal H}(C)$, and then assign standardized
random values to a deterministic orthonormal expansion in $\calh(C)$.
This yields paths as random functions with
an explicit representation formula. 
Using Lo\'eve isometry, we prove that pointwise regularity notions
like continuity or differentiability hold on functions
of $\calh(C)$, paths, and the random field $R_C$ in precisely the same way. Yet, norms of paths as functions behave differently, as we prove that
paths are {\em a.e.} not in $\calh(C)$, but in certain {\em larger}
spaces that can partially be characterized. In case of \M{} kernels
generating Sobolev space $W_2^m(\R^d)$, paths 
lie almost surely in all $W_2^{p}(\R^d)$ for $p<m-d/2$, but almost surely not
in $W_2^{m-d/2}$. This regularity gap between function  and paths
is explained easily by square summability of expansion coefficients
of functions, not of paths. The required orthonormal expansions, well-known
in the probabilistic and the deterministic literature, are analyzed and
compared with respect to convergence rates. 
\end{abstract}


\end{frontmatter}


\section{Introduction}

\subsection{The Problem}
The literature about {\em paths} of random fields is ubiquitous. Reconstruction of the plethora of contributions from mathematics, statistics, machine learning and numerical analysis is beyond the scope of this paper. Nonetheless, this paper is about {\em paths},
{ their definition and their regularity properties.} \\
{To define paths}, it is customary to start with  a {\em stochastic process } or {\em random field} $R$ on a set  $\Omega$, as a mapping that assigns a random variable $R(x)$ to each $x\in \Omega$. The random variables have finite variance, with zero mean,  and we stick to these assumptions throughout the paper. The set $\Omega$
will here be a bounded domain in $\R^d$.  
\biglf
The fact that the random field has finite variance ensures the {\em covariance function}
 $C\;:\;\Omega\times\Omega\to \R$ being equal to
 \begin{equation} \label{eqCxy}
C(x,y):=\text{cov}(R(x),R(y))=\E (R(x)\cdot R(y)),
\end{equation} 
for all $x,y\in \Omega$, to be well defined. Yet, it is hardly ever explained how the expectation $\mathbb E$ above is to be taken - see the attempt in Chapter I of \cite{christakos2013random}. In probability theory, such an expectation is taken over the joint probability space for the pair $(R(x),R(y))$, but in case of 
a random field over $\Omega$ it should be the marginalization of a joint probability space for the full set $\{R(x)\;:\;x\in \Omega\}$. In general, the above covariance can depend on all other $R(z)$, {$z \in \Omega$}, and then (\ref{eqCxy}) does not make sense, notationally.
This is where problems start for the general situation.
\biglf
Things get even worse when {\em paths} are to be defined. By general consent, they are  simultaneous realizations of all $R(x)$ at the same instant, but how does this work? Nature is able to provide a temperature or wind velocity field simultaneously at all points of Earth and all given times, but what is the mathematical technology to simulate that? How can one throw uncountably many strongly dependent dice at the same time? 
\biglf
Here, the general consent is to say that the map 
$x\mapsto R(x,\xi)$ is a {\em path}, where $\xi$ is a real value from the state space $\R$  of the identical probability space of all $R(x)$. This would mean that paths are functions on $\Omega\times \R$ and
conceals that there is a draw from the joint probability space in the background that yields itself a function of $x$. This notational nonsense even made it into the wikipedia \citep{wiki:sto-pro}, following plenty of books and papers that are not cited here. For each fixed realization $\xi$, now from the joint probability space, whatever it may be, the deterministic function 
$x\mapsto R(x,\xi(x))$  should be called a
path $R_\xi$ of $R$. {This should be the base to study, {\em e.g.}, the regularity properties of sample paths. To make an example, } imagine to take derivatives of it with respect to $x$.
\biglf
Because standard probability distributions do not exist on infinite
dimensional spaces  {\citep{oxtoby:1946-1}}, both steps work directly on finite sets $\Omega$ only, and
require a workaround for  infinite sets. This is the main goal of this paper. 
\subsection{The Turnaround}
\noindent {To define paths,} the standard approach performs the steps
\begin{tcolorbox}
\begin{center}
\resizebox{\textwidth}{!}{  
\begin{tabular}{cccccc}
Random Field $R$ &$\Rightarrow$ & Covariance $C$ &$\Rightarrow$ & Paths $p$
\end{tabular}
}
\end{center}
\end{tcolorbox}
\noindent where {the first two} steps require a nonexisting joint distribution. 
In contrast to this, our turnaround is
\begin{tcolorbox}
\begin{center}
\resizebox{\textwidth}{!}{  
\begin{tabular}{cccccc}
Covariance $C$ &$\Rightarrow$ & Paths $p$  &$\Rightarrow$ & Random Field $R$ 
\end{tabular}
} 
\end{center}
\end{tcolorbox}
\noindent with two intermediate steps to be described below. It will avoid joint distributions, and it models Nature better, because it starts from the global picture. It should be the standard modelling technique. Before details of the Turnaround are provided in Section \ref{SecTurn}, the standard workaround should be presented.
\subsection{The Gaussian Workaround}
Gaussian processes {are} 
extremely popular in all disciplines, including approximation theory, statistics, and in particular machine learning {\citep{williams1995gaussian, williams2006gaussian, seeger2004gaussian}}. Their main feature is that the joint distribution of finitely many random variables is multivariate Gaussian, and all marginal and conditional distributions are Gaussian again.
\biglf
This leads to a standard workaround, {exploiting }
the fact that when the Gaussian process is
restricted to a finite point set $X\subset\Omega$,
the covariance function is restricted to $X\times X$
and inserted as a matrix into the multivariate Gaussian distribution.
If $x$ and $y$ are elements of $X$, the covariance of $R(x)$ and $R(y)$ is then
given by $C(x,y)$ and is independent of $X$.
Therefore, one can work bottom-up on large finite sets
in the Gaussian case
{\citep{billingsley1995probability}. The }
smart trick is to define paths through limits on countable dense subsets of $\Omega$, avoiding the notion of 
a joint Gaussian distribution on $\Omega$.  Hence, certain properties {of paths} can be defined through continuation to $\Omega$. The strategy now is 
\begin{tcolorbox}
\begin{center}
\begin{tabular}{ccc}
Covariance $C$ &$\Rightarrow$ & joint distribution \\  
&& \\
& $\Rightarrow$ & Random Field $R$ \\ && \\ &$\Rightarrow$&Paths $p$ 
\end{tabular}
\end{center}
\end{tcolorbox}
\noindent
but restricted to finite sets for the first two steps and a
continuity argument for the last.
\biglf
Now the definition of paths needs a transition from finite
to infinite point sets. 
And since regularity properties are based
on limits, this transition needs special care because it involves a double limit. However, the Gaussian workaround is successful and presented widely in the literature.  
\biglf
In the course of the paper, it will be shown that the Turnaround is fully compatible with the Gaussian work\-around, while simplifying all arguments considerably. In particular, our path definition leads to an explicit formula \eqref{eqpSWR}, and  {Theorem \ref{TheTruncErr}} will handle the path continuation problem. 
\biglf
There is a highly technical and abstract workaround towards abstract Wiener space that we do not follow here. It is well described in \cite{stroock2023gaussian}. 
\subsection{Roadmap}
The next section will provide the details of the proposed Turnaround. In particular, two intervening steps are needed: the Hilbert space $\calh(C)$ where the covariance kernel $C$ is reproducing, and the orthonormal expansions therein, which will be the basis for defining paths first and then the random field $R_C$.
\biglf
Based on this, pointwise regularity theorems for random fields $R_C$ and functions
from the native space $\calh(C)$ are treated 
{simultaneously} in Section \ref{SecPwReg}.
A basic tool is the Lo\'eve isometry from Section \ref{SecLI}
that comfortably relates the deterministic Hilbert space $\calh(C)$
to an isometric Hilbert space $\calh(R_C)$ of second-order zero-mean random
functions. Deterministic and probabilistic results
turn out to be in total agreement, if interpreted appropriately,
but they do not address paths.
\biglf
Then Section \ref{SecPoRF} turns to paths, pointwise, as functions,
and with norms of paths
in function spaces as random variables.
Pointwise regularity of paths follows the previous section and
is proven to agree
with pointwise regularity of functions in $\calh(C)$ and regularity
of the random field $R_C$. 
Yet, a gap appears between
regularity of paths as functions
and functions, generated by the same covariance kernel $C$.
The crucial point is that the coefficients of  paths  \eqref{eqpSWR}
are not square summable almost surely, while the corresponding
expansions for functions in $\calh(C)$ are. Exploiting convergence results
for orthonormal expansions, this allows to characterise
the larger spaces in which the paths live.
\biglf
As an interlude, Section \ref{SecOE} recalls the standard orthonormal
expansions in Hilbert spaces  and of random fields.
They come either eigensystem-based
as Mercer or Karhunen-Lo\'eve expansions, or point-based
as Newton expansions or Cholesky decompositions of kernel matrices.
The former have optimal convergence properties, while the latter depend
on how well points are chosen. But if points are chosen by
the $P$-greedy strategy known from {
\cite{mueller-schaback:2009-1,pazouki-schaback:2011-1}},
the Newton/Cholesky decompositions have an asymptotically equivalent
convergence rate 
{\citep{santin-schaback:2016-1,santin-haasdonk:2017-1}}.
{\section{The Turnaround}\label{SecTurn}}
\biglf 
\noindent {{Our Turnaround} consists of the following steps:}
\biglf 
\begin{description}
    \item[Deterministic Phase:] 
    \begin{enumerate}
        \item[D1] We start with a continuous, symmetric and
  positive definite kernel $C\;:\;\Omega\times \Omega\to \R$. This is done in {Section \ref{SecKer}}; 
  \item[D2] Kernels are reproducing in a Hilbert space of functions on $\Omega$, denoted ${\cal H}(C)$ throughout. This is explained through Section \ref{SecHilSpa};
  \item[D3] We consider orthonormal systems {(ONS)} 
  $W=\{w_n\}_{n\in\N}$ in $\calh(C)$
  in Section \ref{SecOE}. They have an explicit connection with $C$ and ${\cal H}(C)$ in points 1 and 2 through 
\begin{equation} \label{eqSwn2C}
\displaystyle{ \sum_{n\in\N} w_n^2(x)=C(x,x),\;x\in \Omega.  } 
\end{equation}
\end{enumerate}
\item[Transition to Randomness:] The second phase will allow for a proper definition of paths that covers a central role in this manuscript. Specifically, 
\begin{enumerate}
    \item[R1] In Section \ref{SecEPRF}, we
select a univariate
standard probability space over $\R$ with Lebesgue measure. We then take
 independent samples, denoted $S$ throughout, and being sequences $S=\{s_n\}_{n\in\N}$ in
$\R^\N$ where each component is an independent realization of a
random variable following a zero-mean variance-one distribution
$\calr$. 
\item[R2] This allows to define {\em $(S,W,\calr)$-Expansion Paths}, denoted $p_{S,W,\calr}$, and specified through the identity 
\begin{equation} \label{eqpSWR}
p_{S,W,\calr}=\displaystyle{\sum_{n\in\N} s_nw_n.  } 
\end{equation}
\item[R3] This implies, for every $x$, the existence of a random variable $R_{C}(x)$ {whose samples are
\begin{equation} \label{eqRCsam}
R_{C,W,\calr}(x)=p_{S,W,\calr}(x)=
\displaystyle{\sum_{n\in\N} s_nw_n(x),  } 
\end{equation}
and it will turn out that 
$$
\text{cov}(R_{C,W,\calr}(x),R_{C,W,\calr}(y))=C(x,y) 
$$
holds for all $x,y\in \Omega$. }
\end{enumerate}
\end{description}
Note that this proceeds from $C$
to paths first, and then to random fields. {In more detail:}
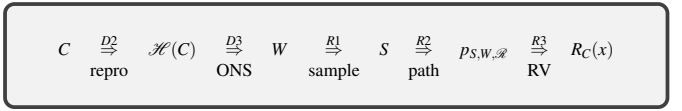
\begin{tcolorbox}
\begin{center}
\resizebox{\textwidth}{!}{  
\begin{tabular}{ccccccccccc}
$C$ 
&$\stackrel{D2}{\Rightarrow}$ & $\calh(C)$
&$\stackrel{D3}{\Rightarrow}$ & $W$
&$\stackrel{R1}{\Rightarrow}$ & $S$ 
&$\stackrel{R2}{\Rightarrow}$ & $p_{S,W,\calr}$ 
&$\stackrel{{R3}}{\Rightarrow}$ & $R_C(x)$ \\
&repro && ONS && sample && path && RV
\end{tabular}
}
\end{center}
\end{tcolorbox}
It works in general because
\eqref{eqSwn2C} implies that $p_{S,W,\calr}(x)$ is finite almost surely. 
{In contrast to the standard approach,
the realizations 
of the random variables  $R_C(x)$ are {\em exactly}  the path values at $x$.} Expansions are used widely in the literature, but not as starting points for defining paths and random variables. 
\biglf
The special case of a Gaussian random field is attained by setting $\calr=\caln(0,1)$ without needing a detour via
limits of point sets. These are hidden in the construction of the orthonormal
system, when applying either the standard Cholesky-based sampling procedure
or a Mercer expansion of the kernel, see Section
\ref{SecOE}. In the latter case, the
result is a Karhunen-Lo\'eve
expansion of the random field.  
\biglf
{Expansions are abundantly used in the literature
when it comes to investigate path properties, e.g. in \cite{steinwart:2019-1}, but they are not used for the definition of paths. When starting from a random field $R$, expansions lead to a random field $\tilde R$ that has to be related to $R$, e.g. by being called a ``modification''. Starting from $C$ and expansions eliminates this problem.  }
\biglf
The turnaround will be shown to work flawlessly, 
lose nothing  compared to the standard workaround via Gaussian Processes, but allowing more general distributions
and an easy access to  regularization results concerning paths. 
\section{Deterministic Phase in the Way to Paths}\label{SecBasics}
As announced above, we proceed here from covariance functions $C$ as kernels
to the Hilbert spaces $\calh(C)$ they generate, including their orthonormal
expansions. In Section \ref{SecEPRF}, randomness comes into play and leads to paths and random fields. 
\subsection{D1: Kernels}\label{SecKer}
Covariance functions are called {\em kernels} in the deterministic literature
{\citep{aronszajn:1950-1,buhmann:2003-1,wendland2004scattered,fasshauer-mccourt:2015-1}} and are real-valued mappings
$$
C\;:\;\Omega\times \Omega\to \R.
$$
They are {\em symmetric}, i.e. $C(x,y)=C(y,x),\;x,y\in \Omega$ and
in most cases {\em positive definite}. This means that for each finite subset
$\Xn$ of $\Omega$ the symmetric $n\times n$ {\em kernel matrix} $\boC_{X_n,X_n}$
with entries $C(x_j,x_k),\;1\leq j,k\leq n$ 
is positive definite. 
A kernel $C$ on $\Omega$ is called {\em stationary}
or {\em translation-invariant} 
if it is a function of the  {\em lag}  $x-y$.
We shall be sloppy for writing $C$ again for a
kernel that is translation invariant, and we always assume positive
definiteness.  
\biglf 
A very common additionally assumption for stationary kernels is 
{\em radial symmetry} or {\em isotropy}:
\begin{equation}
\label{covariance}
C({x}-{y})= \sigma^2 \varphi(\|{x}-{y}\|),
\end{equation}
for ${x},{y} \in \R^d$ and $\|\cdot\|_2$
denoting the Euclidean distance.
For convenience, the function
$\varphi$ is continuous and normalized, so that $\varphi(0)=1$.
The parameter $\sigma^2$ 
then is the constant variance of $R$. In many cases, the scalar
$r=\|x-y\|$ is used for the shortcut
$C({r})= \sigma^2 \varphi(r)$. Whenever no
confusion can arise, the word isotropy will be omitted whenever
it is apparent from
the context. We use the words
{\em covariance function} or {\em covariance kernel} if there is a
probabilistic background, and  {\em kernel} elsewhere,
keeping the notation $C$. 
\biglf 
We now describe the r\^ole of kernels as covariance functions,
{ignoring their importance for machine learning
\citep{schoelkopf-smola:2002-1,sutton2012introduction}
and meshless methods for solving partial differential equations
\citep{schaback-wendland:2006-1}}.
\biglf
Any isotropic covariance in $\R^d$ has a representation as a nonnegative mixture of the type \citep[Theorem 1]{schoenberg:1938-1}:
\begin{equation*}
    \varphi(r) =  2^{\frac{d}{2}-1} \Gamma\left(\frac{d}{2}\right) r^{1-\frac{d}{2}} \int_0^\infty z^{1-\frac{d}{2}} J_{\frac{d}{2}-1}(zr) \text{d}F_d(z), 
\end{equation*}
$r>0$, where $F_d$ is a nondecreasing bounded measure on $(0,+\infty)$. We follow \cite{daley-porcu} to call $F_d$ a $d$-Schoenberg measure throughout. When $\varphi$ is absolutely integrable, then $F_d$ is absolutely continuous with respect to the Lebesgue measure, with radial spectral density $\widehat{\varphi}$ that is attained through
\begin{equation*}
    \widehat{\varphi}(z) = \frac{1}{(2\pi)^{\frac{d}{2}}} z^{1-\frac{d}{2}} \int_0^{+\infty} r^{\frac{d}{2}} J_{\frac{d}{2}-1}(z r) \varphi(r) {\rm d}r, \qquad z > 0.
\end{equation*}

An very important special case of random
fields has isotropic \M{} covariance functions 
\begin{equation}\label{matern0}
  {\cal M}_{\nu, \alpha} (r) =
  \frac{2^{1-\nu}}{\Gamma(\nu)} \left ( \frac{r}{\alpha} \right )^{\nu} {\cal
    K}_{\nu}
  \left ( \frac{r}{\alpha} \right ), \qquad r \ge 0,
\end{equation}
where $\alpha > 0$ is a {scale} parameter, $\nu >0$ is
called the \emph{smoothness} parameter,
and ${\cal K}_{\nu}$ is a modified Bessel function of the second kind of order $\nu$. 
The isotropic {Mat{\'e}rn}  spectral density or radial Fourier transform is
\citep[11.4.44]{abramowitz-stegun:1970-1}
\begin{equation} \label{stein0}
\widehat{{\cal M}}_{\nu,\alpha}(z)= \frac{\Gamma(\nu+d/2)}{\pi^{d/2} \Gamma(\nu)}
\frac{ \alpha^d}{(1+\alpha^2z^2)^{\nu+d/2}}
, \; \; z \ge 0, 
\end{equation}
where $z=\|\omega\|$
is the scalar variable in Fourier space corresponding to $r$.
See \cite{porcu2023mat} for an overview of 
random fields with \M{} covariance functions, among other cases with similar spectral behaviour.
\biglf
The {\em Generalized Wendland}
\citep{Gneiting:2002,zastavnyi2006some} kernel is defined for $\kappa, \beta>0$, as 
\begin{equation} \label{WG2*}
{\cal W}_{\mu,\kappa,\beta}(r):= 
\frac{1}{B(2\kappa,\mu+1)} \int_{r/{\beta}}^{1} u(u^2-(r/\beta)^2)^{\kappa-1} (1-u)^{\mu}\,{\rm d} u  , 
\end{equation} when $ r/\beta < 1$, and $0$ otherwise. 
with $B$ denoting the beta function. The univariate kernel  ${\cal W}_{\mu,\kappa,\beta}(r)$ is positive definite in $\mathbb{R}^d$ 
if and only if $\mu \ge (d+1)/2+ \kappa$ \citep{zastavnyi2006some}. The special case $\kappa=0$ is known as
the Askey family \citep{Askey:1973} 
\begin{equation*} 
{\cal A}_{\mu,\beta}(r) :=   \begin{cases}   \left ( 1- r/\beta \right )^{\mu} ,& 0 \leq r/\beta< 1,\\ 0,&r/\beta \geq 1. \end{cases}
\end{equation*}
We define ${\cal W}_{\mu,0,\beta}:= {\cal A}_{\mu,\beta}$. 
\biglf
Algebraically closed
form solutions for \eqref{WG2*} are available when
$\kappa=k$, a positive integer.  We have
 $ {\cal W}_{\mu,k,\beta}(r) =
{\cal A}_{\mu+k,\beta}(r) P_{k}(r)$,
with $P_{k}$ a polynomial of order $k$.
Such a case is termed {\em original} Wendland functions after the tour de force
by   \cite{Wendland:1995}. {Then }
\cite{Schaback:2011} proved that other closed form
solutions can be obtained when $\kappa=k+0.5$. {These were called }
\emph{missing Wendland}  functions. Arguments in Theorem
1(3) of \cite{bevilacqua2019estimation} prove that the Fourier transform of the
generalized Wendland function, $\widehat{{\cal W}}_{\mu,\kappa,\beta}$
is not as simple as (\ref{stein0}), but has a similar 
asymptotic order of $z^{-2 \lambda}$ for $z\to\infty$
provided $\mu \ge \lambda:= (d+1)/2+\kappa$. 
\biglf
Mat{\'e}rn and Generalized Wendland spectral
densities can be parameterized in such a
way that their tails behave similarly.
Specifically, Theorem 1 in \cite{bevilacqua2019estimation} shows that $2 \nu = 2
\kappa +1$ is needed. The scale factors $\alpha$ and
$\beta$ in the two models
have no influence on the spectral asymptotics.
Yet, they play a key role, in concert with
smoothness and variance, to determine conditions for equivalence of Gaussian
measures associated to random fields with
these classes of covariance functions
\citep{bevilacqua2019estimation}.
In turn, such conditions are the crux
to provide sufficient condition for best
optimal linear prediction under a
misspecified covariance function \citep{stein:1999-1}.
\subsection{D2: Hilbert Spaces}\label{SecHilSpa}
{We now come back to general kernels and state their main properties
  as needed for the paper.} 
The space $\calh(C)$ spanned by all {\em kernel translates}
$C(x,\cdot),\;x\in\Omega$ is a Hilbert space of functions
on $\Omega$ that has an inner
product $(.,.)_{\calh(C)}$ with the properties 
\begin{equation} \label{eqHCrepro}
\begin{array}{rcll}
f(x) &=& (f, C(x,\cdot))_{\calh(C)},&f\in{\calh(C)},\;x\in\Omega,\\
&& \\
C(x,y)&=& (C(x,\cdot),C(y,\cdot))_{\calh(C)},&x,y\in\Omega,\\
&& \\
\lambda^x\mu^y C(x,y) &=& (\lambda,\mu)_{{\calh(C)}^*},&\lambda,\mu\in{\calh(C)}^*\\
\end{array}
\end{equation}
where ${\calh(C)}^*$ denotes the (topological)
dual of ${\calh(C)}$ containing all continuous (i.e. bounded) linear functionals
on ${\calh(C)}$ and $\lambda^x$ means evaluation
of $\lambda$ with respect to the variable $x$.
In particular, all point evaluation functionals $\delta_x\;:\;f\mapsto f(x)$ are
continuous with norm $\sqrt{C(x,x)}$. {Details are in monographs by  \cite{schaback:1997-3,abdelaziz-hamouine:2008-1,wendland2004scattered,fasshauer-mccourt:2015-1}}.
\biglf
For Sobolev spaces, we work on $\R^d$ and use the definition
$$
W_2^m(\R^d):=\displaystyle{\left\{f\;:\;\R^d\to \R, \int_{\R^d}|\hat f(\omega)|^2(1+\|\omega\|_2)^m d\omega<\infty   \right\}}
$$
as subspaces of $L_2(\R^d)=W_2^0(\R^d)$. 
They are separable Hilbert spaces under the inner product
$$
(f,g)_{W_2^m(\R^d)}=\displaystyle{\int_{\R^d}\hat f(\omega)\overline{\hat g(\omega)}(1+\|\omega\|_2)^m d\omega}, 
$$ for $f,g\in W_2^m(\R^d)$,
and there is a reproduction formula
$$
f(x)=(f,{\cal M}_{m-d/2,1}(\|x-\cdot\|_2)_{W_2^m(\R^d)},
$$ for $f\in W_2^m(\R^d),\,x\in \R^d$, by
using a special version of the general isotropic \M{} kernel
from (\ref{matern0}). {This links the most important function spaces
  of Real Analysis and Partial Differential Equations to the most important
  covariance functions in Spatial Statistics.}
\biglf
{In turn, Theorem 5 in \cite{bevilacqua2019estimation} ensures that the above reproducing property is guaranteed provided $m-d/2 = \kappa + 1/2 $, that is when $\kappa = m - (d-1)/2$. }
\biglf
Summarizing, each positive definite
{covariance function} $C$ is a
{kernel} that leads to a {\em native} Hilbert space ${\calh(C)}$ of
deterministic functions on $\Omega$, and Sobolev spaces arise from {
  \M{} kernels}. 
\subsection{D3: Orthonormal Expansions}\label{SecOE}
{Like in all separable Hilbert spaces, there are orthonormal systems
$\{w_n\}_{n\in\N}$ in $\calh(C)$ that allow each function $f\in\calh(C)$
to be written as a series
$$
f=\displaystyle{\sum_{n\in \N}(f,w_n)_{\calh(C)}w_n    } 
$$
with
$$
\|f\|^2_{\calh(C)}=\displaystyle{\sum_{n\in \N}(f,w_n)^2_{\calh(C)}   } 
$$
and convergence in norm. In particular,
$$
C(x,y)=\displaystyle{\sum_{n\in \N}w_n(x)w_n(y)    } 
$$
holds for all $x,y\in \Omega$, and there is a pointwise bound
\begin{equation} \label{eqwnC}
\displaystyle{\sum_{n\in \N}w_n(x)^2 =C(x,x)   } \fa x\in \Omega 
\end{equation}
for the functions $w_n$ that ensures $w_n(x)\to 0$ for $n\to\infty$
and
$$
\displaystyle{\sum_{n\in \N}\|w_n\|^2_{L_2(\Omega)} =\int_\Omega C(x,x)dx   }. 
$$
We now devote some attention to two special cases of orthonormal expansions. 
\biglf
The standard way of generating path approximations for random fields
is by a stepwise Cholesky decomposition of kernel matrices.
We show here that it yields a special case of an orthonormal expansion.
\biglf
Consider a countable dense set  $X_\infty$ of points in $\Omega$.
Then the Cholesky decomposition procedure for
infinite kernel matrices $\boC_{X_\infty,X_\infty}$ can be written as 
$$
C(x_j,x_k)=
\sum_{n=1}^\infty w_{n,j}w_{n,k},\;1\leq j,k <\infty 
$$
with an infinite triangular matrix $\boW$ with entries $w_{n,j}$ that satisfy $w_{n,j}=0,\;1\leq j\leq n-1$ and $w_{nn}=:\sigma_n>0,\;n\geq 1$. The standard stepwise Cholesky procedure generates exactly those matrices for $n\to\infty$, and the result is $\boC_{X_\infty,X_\infty}=\boW^T\boW$.
\biglf
From a more general viewpoint, this is an expansion
\begin{equation}\label{ExpGen}
 C(x,y) =\displaystyle{\sum_{n=1}^\infty w_n(x)w_n(y)}  \fa x,y\in\R^d
\end{equation}
of the kernel itself into functions $w_n,\,n\geq 1$.
The functions $w_n$ have the properties 
$$
\begin{array}{rcll}
  w_n(x_j) &=& 0,&1\leq j<n\\
  w_n(x_n) &=:& \sigma_n>0,\\
  w_n & \in & span\{C(x_1,\cdot),\ldots,C(x_n,\cdot)\},& n\geq 1.
\end{array} 
$$
In the deterministic literature 
\citep{mueller-schaback:2009-1,pazouki-schaback:2011-1}
this is the well-known expansion into a {\em Newton basis}, and it is orthonormal in $\calh(C)$. The recursion formulas
in terms of full functions are
\begin{equation} \label{eqwrec}
\begin{array}{rcl}
w_j^2(x_j)&=&C(x_j,x_j)-\displaystyle{\sum_{m=1}^{j-1}w_m(x_j^2)   },\\ 
w_j(x)w_j(x_j)&=&C(x_j,x)-\displaystyle{\sum_{m=1}^{j-1}w_m(x_j)w_m(x)   }.\\ 
\end{array} 
\end{equation} 
Note that expansions  (\ref{ExpGen}) replace a stationary kernel by a sum of
products, and stationarity of is lost for partial sums.
\biglf
This ends the deterministic phase. We now have all tools
needed for paths \eqref{eqpSWR} and random variables \eqref{eqRCsam}.

\section{The Stochastic Phase: 
Expansion Paths and Random Fields}\label{SecEPRF}
After starting from kernels $C$ and providing their deterministic
properties, it is now time to introduce randomness. This will turn the kernel
into a covariance function for a random field,
without using a joint distribution. We shall use the orthonormal expansions of Section \ref{SecOE} for this purpose.
\biglf
{Step R1 introduces } a zero mean and unit variance distribution, ${\cal R}$, from which an independent and identically distributed sequence, $S=\{s_n\}_n$ is drawn. The background is a completely different random field, now on $\N$ and denoted by $\N^\calr$, with a trivial joint distribution induced by independent ${\cal R}$-distributed random variables $S(n)$ for all $n\in \N$.   
Its paths $S=\{s_n\}_n$ are well-defined due to the independence, and are the sequences $S$ used here. 
\biglf
{Step R2 {couples} the sample $S=\{s_n\}_n$ {with} an orthonormal expansion $W$ in $\calh(C)$ to arrive at 
\eqref{eqpSWR}.
There, we called the function}  $p_{S,W,\calr}$ an {\em $(S,W,\calr)$-expansion path}. 
\biglf
{Finally, step R3 looks locally at \eqref{eqpSWR} on }
 $x$, to get a random variable $R_{C,W,\calr}(x)$
whose realizations are in \eqref{eqRCsam}
with realizations $S$ from $\N^\calr$. Note that these random variables cannot
be sampled
individually. A sample $S$ from $\N^\calr$ provides samples for all
$R_{C,W,\calr}(x)$ simultaneously. 
As long as we keep $C,W,$ and $\calr$ fixed, we shall reduce
the notation to $R_C(x)$ for reasons to be apparent below.
{
\subsection{Properties of Paths and Random Variables}
}
Because we first defined paths and then pointwise random variables,
both $R_C(x)$ and $R_C(y)$ have to use the same samples $S$, and then 
\begin{equation} \label{eqcRRC}
\begin{array}{rcl}
{\rm cov}(R_C(x),R_C(y))&=&\E\left(\displaystyle{\sum_{n\in \N}s_nw_n(x)}\right)
\left(\displaystyle{\sum_{m\in \N}s_mw_m(y)}\right)\\
&=&
\displaystyle{\sum_{n\in \N}w_n(x)w_n(y)}\\
&=&
C(x,y)
\end{array} 
\end{equation}
for all $x,y\in \Omega$. 
\biglf
{The above approach makes it easy to construct non-Gaussian
  random fields with a prescribed covariance function, $C$. Users still have the
  orthonormal system $W$ and the standardized probability distribution $\calr$
  at their disposal. 
  \biglf
  It is interesting to study the set of all
  paths \eqref{eqpSWR} for fixed $C,\,W$, and $\calr$. By the following theorem,
  it is independent
  of $W$ as a set, but different $W$ and $\calr$ will prefer
  certain subsets of paths over others when producing repeated paths.
  \begin{Theorem}\label{TheWIndep}
    The set of all possible expansion paths \eqref{eqpSWR} is independent of the
    orthonormal expansion $W$, if $C$ and {$\calr$} are fixed.
\end{Theorem} 
\begin{proof}
If there are two orthonormal systems $W$ and $U$ with
$$
w_n=\sum_{m\in \N} d_{nm}u_m,\;n\in \N, 
$$
the infinite transformation matrix $D$ is orthogonal, i.e.
$$
DD^T=D^TD=Id,
$$
letting the series be absolutely convergent.  Now 
path representations  $p_{S,W,\calr}$ and  $p_{T,U,\calr}$ can be 
related by  
$$
t_m=\displaystyle{\sum_{n\in\N}s_n d_{nm}   }, \;m\in\N
$$
to show that
$$
\begin{array}{rcl}
  p_{S,W,\calr}
  &=&
  \displaystyle{\sum_{n\in \N}s_nw_n    }\\
  &=&
  \displaystyle{\sum_{n\in \N}s_n\sum_{m\in \N} d_{nm}u_m   }\\
  &=&
  \displaystyle{\sum_{m\in \N}u_m\sum_{n\in \N}s_n d_{nm}  }\\
  &=&
  p_{T,U,\calr}.
\end{array} 
$$
\end{proof} 
\subsection{Statistical Equivalence}\label{SecStatEv}
Assume two orthonormal systems $v=\{v_n\}$ and $w=\{w_n\}$ and keep the same
univariate distribution $\calr$ and the same
covariance function $C$.
The two systems could be called {\em statistically
  equivalent}, if for all
measurable sets $M\subset \Omega\times \R$ of finite Lebesgue measure,
the probability of paths based on $v$ and $w$  to have their graphs in $M$
should be the same. 
\biglf
A simplification replaces $M$ by {\em cylinder sets}
$$
\calc(X_n,\boa,\bob):=
\{f\;:\;\Omega\to \R\;:\; a_j\leq f(x_j)\leq b_j,\;1\leq j\leq n\}
$$
for $\boa<\bob\in \R^n$, and this defines a probability 
$$
\text{Prob}\displaystyle{   \left(  a_j\leq \sum_{m\in \N}s_mw_m(x_j)\leq b_j,\;1\leq j\leq n \right)}
$$
that a path based on $w$ hits $\calc(X_n,\boa,\bob)$. This is
\begin{equation} \label{eqpcyl}
\text{Prob} \displaystyle{   \left(  a_j\leq R_C(x_j) \leq b_j,\;1\leq j\leq n \right)}
\end{equation}
and { is independent of  the expansion 
whenever the joint distribution
of all $R_C(x_j), 1\leq j\leq n$ is.
\begin{Theorem}
Under Gaussianity,  
there is statistical equivalence of all expansion paths on all cylinder sets.\qed
\end{Theorem}} 
Probabilities \eqref{eqpcyl} on cylinder sets are well-defined in case of
Gaussianity, and the probability
stays the same if $X_n$ is contained in a point-based expansion
for an infinite dense set $X_\infty$. This why users
can consider the standard constructions of Gaussian paths using point sets $X_\infty$ to be reliable and independent of expansions.
Cylinder sets arise in \cite{billingsley1995probability} and various approaches to abstract Wiener spaces, used to define so-called {\em cylinder set measures} \citep{enwiki:1240567585} which usually are no measures. 
{
\subsection{Paths for Point-Based Expansions}}
If the standard
Cholesky decomposition of the kernel matrix with entries
$C(x_j,x_k),\;1\leq j,k\leq N$ is executed only for a finite point set $X_N$,
the path approximations are the partial sums of \eqref{eqpSWR}treated in
Theorem \ref{TheTruncErr} {below},  but usually only calculated on $X_N$.
The recursion formula \eqref{eqwrec} extends the result to be a function
on all of $\Omega$, and it can easily be proven that
the extension is the Kriging interpolant using the
finite path vector on $X_N$.
Theorem \ref{TheTruncErr} {will show} how far this is
from a calculation on $\Omega$ itself. The {\em residual} kernel
\begin{equation} \label{eqCCww}
C_{N+1}(x,y):=C(x,y)-\displaystyle{\sum_{n=1}^Nw_n(x)w_n(y)   } 
\end{equation}
is the kernel $C$ conditioned to $X_N$, and $C_{N+1}(x,x)$ is the variance for
Kriging on $x$ given values on $X_N$,  {because it coincides with the
{\em Power Function} that has this property,  \citep[see, {\em e.g.},][p. 98]{fasshauer-mccourt:2015-1}. The early paper
\citep{wu-schaback:1993-1} called it {\em Kriging function}}. This is the probabilistic interpretation
of such a point-based expansion. The presentation here is in terms of functions,
not finite path vectors. The limit $N\to\infty$ causes no problems
if points are reasonably selected to become dense in the limit. \citep{mueller-schaback:2009-1}.
\\
\subsection{Paths for Mercer Expansions}
For a Mercer {\citep{mercer:1909-1}}
expansion {of the kernel}, the basic equations are
\begin{equation}\label{eqMercer}
\int_\Omega C(x,y)v_n(x)dx=\lambda_n v_n(x),\; x\in \Omega
\end{equation}
with positive $\lambda_n$ decaying to zero and
orthonormality of the $v_n$ in $L_2(\Omega)$ paired
with orthogonality $(v_n,v_m)_{\calh(C)} =\delta_{nm}\lambda_n^{-1},\,n,m\geq 1$
in the native space $\calh(C)$. Then
$$
C(x,y)=\sum_{n=1}^\infty \lambda_n v_n(x)v_n(y)
$$
is an expansion in $L_2(\Omega)$. We get (\ref{ExpGen}) and orthonormality in $\calh(C)$  by setting $w_n=\sqrt{\lambda_n}v_n$.
In the probabilistic context, this turns into a  Karhunen-Lo\`eve expansion.
Paths have the form \eqref{eqpSWR}, in particular
$$
p_{S,W,\calr}=\displaystyle{\sum_{n\in \N}s_nw_n = \sum_{n\in
    \N}s_n\sqrt{\lambda_n} v_n   }
$$
with norm
$$
\|p_{S,W,\calr}\|_{L_2(\Omega)}^2=\displaystyle{\sum_{n\in \N}s^2_n\lambda_n   }.
$$
A \MKL{} expansion
is not point-based, and therefore it lacks the above probabilistic
interpretation. But it is optimal in the sense that the decay of the
$\lambda_n$ cannot be improved given $C$ and the space $\calh(C)$,
leading to optimal convergence of \MKL{} paths
{\citep{santin-schaback:2015-1}}.
\biglf
If points in a point-based expansion are chosen by
the {\em $P$-greedy method}
{\citep{mueller-schaback:2009-1}}, picking $x_{N+1}$ as
$$
x_{N+1}=\arg \max \{ C_{N+1}(x,x),\;x\in\Omega\},
$$
the convergence of path approximations is asymptotically optimal as well {\citep{santin-haasdonk:2017-1}}.
This is computationally much cheaper than a
\MKL{} expansion.


\section{{Regularity of Random Fields after the Turnaround}}
\label{SecPwReg}

\subsection{A tool: the Lo\`eve Isometry}\label{SecLI}
{By \eqref{eqcRRC}, the linear space 
\begin{equation} \label{eqSR}
{{\cals(R_C)}}=span\{R_C(x)\;:\;x\in\Omega   \}. 
\end{equation}
of zero-mean random variables
is isometrically isomorphic to the span of kernel translates
$C(x,\cdot)$ by the {\em Lo\'eve} map
$$
\call_{C}(R_C(x)):=C(x,\cdot) \fa x\in \Omega
$$
on the generators. We denote the completion of $\cals(R_C)$
by $\calh(R_C)$ to distinguish it from its
isometric counterpart $\calh(C)$.
If one does not start from $C$, and without Gaussianity,
the definition of the Lo\'eve map
is problematic due to the definition of the covariance
via marginalization of the joint distribution of all $R(\cdot)$ to  
the joint distribution of pairs $(R(x),R(y))$.  
\biglf
We know the Hilbert space completion $\calh(C)$
of the translates $C(x,\cdot),\;x\in \Omega$, and then each
deterministic function $f\in \calh(C)$ defines a zero-mean
second-order random variable $S_f=\call^{-1}_C\in \calh(R_C)$.
It has the property
$$
\text{cov}(S_f,R_C(x))=f(x),\;x\in \Omega
$$
{that can be called the {\em random reproduction formula} that complements the standard deterministic reproduction formula from \eqref{eqHCrepro}, obtained by applying the Lo\'eve isometry.}
\biglf 
But more important are the duals of the two Hilbert spaces. }
Let $\calh(C)^*$ and ${\calh(R_C)}^*$ denote the
dual spaces associated with, respectively,
$\calh(C)$ and ${\calh(R_C)}$. In view of {the}
Lo{\`e}ve isometry, we have that 
functionals on the native space $\calh(C)$
should now correspond to functionals on ${\calh(R_C)}$.
If we start from $\lambda\in\calh(C)^*$,
we can define a functional $\lambda^*\in {\calh(R_C)}^*$ on ${\calh(R_C)}$ by
$$
\lambda(\call_C(S))=:\lambda^*(S) \fa S\in {\calh(R_C)}.
$$
This is the standard duality map $\call_{C}^*$
from $\calh(C)^*$ to ${\calh(R_C)}^*$.
By the Riesz representer theorem, there is a
second-order zero-mean random variable $S_\lambda$ with
$$
\E(S_\lambda S)=\lambda(\call_C(S))=:\lambda^*(S) \fa S\in {{\cals(R_C)}},
$$
and in particular
\begin{equation} \label{eqESR}
\begin{array}{rcll}
\E(S_\lambda R_C(x))&=&\lambda(C(x,\cdot)) \;\;\fa x\in\Omega,\\
\E(S_\lambda S_\mu)&=&(\lambda,\mu)_{\calh^*(C)}
=\lambda^x \mu^yC(x,y)
\end{array}
\end{equation}
for all $\lambda,\,\mu \in\calh^*(C)$.
Here, the superscript $x$ on a functional $\lambda$
denotes action with respect to the variable $x$.
\biglf
The consequence is that
all functionals $\lambda\in\calh^*(C)$ 
on $\calh(C)$ lead to
valid second-order mean-zero random variables $S_\lambda\in {{\calh(R_C)}}$.
Under Gaussianity, all second-order random variables $S_\lambda$ of this type
will be Gaussian again.
{And the admissible and bounded functionals on both Hilbert spaces
  are comparable
  via the dual of the Lo\'eve map.}
This argument will be crucial when investigating
pointwise regularity notions for both {functions in $\calh(C)$ and
  random fields $R_C$ }  in the next section.

{
  In particular, a classical
  pointwise derivative is a functional $\lambda$ on $\calh(C)$, and
  then $S_\lambda$ is a random variable with the same norm that
  describes the corresponding derivative of the random field $R_C$
  in the mean-square sense. We shall describe this in more detail.}
\subsection{Pointwise differentiability} 
Pointwise derivative
functionals $\delta^\alpha_x(f):=
{(D^\alpha f)(x)}$ that
are admissible in the native space $\calh(C)$ are those having a finite norm $\|\delta ^\alpha_x\|^2_{\calh^*(C)}$ that is defined through 
\begin{equation}\label{deltavar}
\|\delta ^\alpha_x\|^2_{\calh^*(C)}=
\delta ^{\alpha,u}_x\delta ^{\alpha,v}_x C(u,v)<\infty,
\end{equation}
where the upper indices $u$ and $v$ indicate the variables the functionals use. This induces derivative kernels 
$$
C^{\alpha,\alpha}(x,y)\;: \;(x,y)\mapsto \delta ^{\alpha,u}_x\delta ^{\alpha,v}_y C(u,v) \fa x,y\in \Omega
$$
that are symmetric and positive semidefinite. Boundedness of pointwise derivative functionals is reduced to existence of ``twin'' derivatives of the kernel. 
\biglf
For  each admissible derivative order $\alpha\in \N^d$ there  is a second-order zero-mean  random variable $S_{\delta ^\alpha_x}$ with the variance in (\ref{deltavar}), and the map $x\mapsto S_{\delta ^\alpha_x} $ is a zero-mean second order random field with covariance function $C^{\alpha,\alpha}$.
Starting from $x\mapsto S_{\delta_x}=R(x)$ for $\delta_x=\delta_x^{(0)}$ we can approximate all higher pointwise derivatives by linear combinations of point evaluations, and then 
the random field  $x\mapsto S_{\delta ^\alpha_x}=:D^\alpha_x (R)$  coincides with the pointwise mean-square 
derivative of the random field $R$.
In this sense we have 
\begin{Theorem}\label{TheRegRF}
  Pointwise mean-square differentiability properties of {the random
    field $R_C$} 
coincide with pointwise differentiability properties
of functions in the native space {$\calh(R_C)$}
of the covariance kernel $C$.\qed
\end{Theorem}
\noindent
{A simple illustration follows. We take $d=1$ for sake of simplicity.
A random variable $R_x'$ is the mean-square
derivative of $R$ at $x$, if 
$$
\E(R_x'-(R(x+h)-R(x))/h)^2\to 0 \hbox{ for } h\to 0.
$$
The above approach defines $R_x'$ via 
$\E(R_x' R(z))=D^{1,u}_xC(u,z)$, $\fa z$ and
$E(R_x' R_x')=D^{1,u}_xD^{1,v}_xC(u,v)$.
Then the above expression is
$$
\begin{array}{rcl}
&&\E(R_x'-(R(x+h)-R(x))/h)^2\\
&=& 
\|\delta_x^1 -\frac{1}{h}(\delta_{x+h}-\delta_x)\|^2_{\calh(C)^*}
\end{array}
$$
and this converges to zero if and only if 
$\frac{1}{h}(\delta_{x+h}-\delta_x)$ converges to
$\delta_x'$ in norm in $\calh(C)^*$.
} 
\subsection{Continuity}
Like continuity of functions, this is not a notion that works with a single
point or a single functional. One way to define
it for a random field $R$ is to say that $R$ is {\em pointwise mean-square
continuous} at $x$ if 
\begin{equation} \label{eqRRxhRx}
\lim_{h\to 0\in\R^d}\E(R(x+h)-R(x))^2=0, 
\end{equation}
and because this is \citep{stein:1999-1}
\begin{equation} \label{eqCC2C}
\lim_{h\to 0\in\R^d} \left (C(x+h,x+h)+C(x.x)-2C(x+h,x)\right)=0,
\end{equation}
it is satisfied if $C$ is continuous.  
When combined with the previous section, we get
\begin{Theorem}\label{TheContRF}
  {Assume a kernel $C$ 
    has continuous derivatives 
  $C^{\alpha,\alpha}(x,y)=D^{\alpha,u}_xD^{\alpha,v}_yC(u,v)$.
  Then the second-order zero mean
  random field
  {$D^\alpha R_C$} is pointwise mean-square
  continuous
  almost surely and has the above covariance function.}
  \qed
\end{Theorem}
Again, the slight 
difference between regularity of random fields and 
functions from their their native spaces lies in the difference
of the two regularity notions, here for continuity.
\section{{Pointwise Convergence of Paths}}\label{SecPoRF}
We now consider pointwise regularity of paths, using our explicit representation (\ref{eqpSWR}).
{We have by definition} 
 \begin{Theorem}\label{TheASPWconv}
  The value of expansion
  paths {$p_{S,W,\calr}$} at a point $x$ is a random variable {$R_C(x)$}
  that has mean zero and variance $\sum_{n=1}^\infty w_n(x)^2=C(x,x)$. {It is finite almost surely.}
\end{Theorem} 
{The proof uses the fact that series of second-order zero-mean random variables converge almost surely, if the variances are summable \citep[Theorem 2.5.6, page 84]{durrett:2019-1}.}
\biglf
A {simple} argument allows us to control the error that is
 committed when taking only partial sums in (\ref{eqpSWR}) up to some $N$.
\begin{Theorem}\label{TheTruncErr}
Path approximations {$p_{N,S,W,\calr}$ by truncated sums
$$
p_{N,S,W,\calr}(x)=\displaystyle{\sum_{n=1}^N w_n(x) {s_n}}
$$
} always lie in the native Hilbert
space generated by the covariance kernel,
and the error {$p_{S,W,\calr}(x)-p_{N,S,W,\calr}(x)$}
is a random variable with variance 
\begin{equation}\label{eqresvar}
C(x,x)-\sum_{n=1}^{N}w_n(x)^2=\sum_{n=N+1}^\infty w_n(x)^2
\end{equation}
for all $x\in\Omega$.
Going from $N$ to $N+1$ decreases the variance  by $w_{N+1}(x)^2$ for all $x$.
\qed
\end{Theorem}
{In a nutshell, this is what allows to go over to infinite point sets
  in the Gaussian case to get paths on all of $\Omega$. Limits of point sets
  as a bottom-up or local-to-global process are replaced here by
  series truncations as a top-down or global-to-local process. Nothing is lost,
  because
  the standard finite-to-infinite path construction for Gaussian processes
goes via an orthonormal system anyway, by Section \ref{SecOE}.}
\biglf
The random variable {$D^\alpha(p_{S,W,\calr})(x)$ on expansion paths
$p_{S,W,\calr}$ from \eqref{eqpSWR} has the representation
$$
D^\alpha(p_{S,W,\calr})(x)=\displaystyle{\sum_{n=1}^\infty D^\alpha (w_n)(x) {s_n}}.
$$
It is a random variable with variance zero and variance } 
$$
\sum_{n=1}^\infty (D^\alpha (w_n)(x))^2=
D^{\alpha,\alpha}C(x,x)=\|\delta^\alpha\|^2_{\calh(C)^*} \fa x\in \Omega.
$$
Therefore Theorem \ref{TheRegRF} on almost sure mean-square
pointwise differentiability also holds for {pointwise values of}
expansion paths, { in the sense that pointwise derivatives 
have finite values almost surely.}
\biglf 
Similarly, we examine pointwise
mean-square continuity via \eqref{eqRRxhRx} and get the condition
$$
\lim_{h\to 0\in\R^d}\E\left(
\sum_{n=1}^\infty {{s_n}}(w_n(x+h)-w_n(x))\right)^2=0, 
$$
which coincides with \eqref{eqCC2C}. Then
Theorem \ref{TheContRF} on almost sure pointwise mean-square
continuity extends to expansion paths. 
\begin{Theorem}
Pointwise regularity and smoothness notions for paths of a random field and functions in the native space of its covariance function coincide. The former are to be understood as random variables with bounded variance. \qed
\end{Theorem}
{
\section{Norm Convergence of 
Paths}\label{SecNorConPath}}
We recall Section {\ref{SecEPRF} for $(S,W,\calr)$-expansion paths
\eqref{eqpSWR}. Once sampled, these are functions like any other deterministic
function, and we can investigate their
regularity with deterministic techniques.}     
\biglf
Each partial sum {of \eqref{eqpSWR}}
is a function in the native space whose squared norm is the sum of squares of the coefficients. For $n \to \infty$, the problem requires instruments from probability theory. When $S$ in \eqref{eqpSWR} is an independent sequence, then $S^2:= \{S^2_n\}_n$ is an independent sequence as well. Hence, we can invoke the celebrated Kolmogorov's {\em Three series  theorem}, for which the series $\sum_n s_n^2$ converges if and only if 
$$ 
\text{Prob} \left ( S_n^2 > K \right ) < \infty \qquad \text{and} \qquad \sum_n \mathbb{E} S_n^2 < \infty,     
$$
where the third condition is not needed in our case because we are working with zero-mean random variables. This provides an immediate implication for the case of {IID} (independently and identically distributed) sequences, which is largely used in machine learning as explained by \cite{scholkopf2022causality}. Although the implication is straightforward, we formalize it below for the convenience of the reader.

\begin{Theorem}\label{Thegap1}
Let the sequence $S$ in \eqref{eqpSWR} be {\em IID}. Then, expansion paths \eqref{eqpSWR} do not lie, almost surely, in the native Hilbert space of the
  covariance kernel.
\end{Theorem}
Some commments are in order. The assumption of ID can be relaxed to provide situations where the expansion paths \eqref{eqpSWR} lie, almost surely, in their Native space. Yet, these are more mathematical artifacts rather than real situations. For example, one might set $S_n^2 \sim {\cal N}(0,1/n^2)$ and retrieve $S_n$ by backwards transformation. Yet, this has no sense for practical applications. Hence, the general message is that expansion paths based on independent sequences fall almost surely outside their original Native space. \\ 
We note that the assumption of independence as per Kolmogorov can be relaxed at the price of very technical conditions \citep{brown1971general}. Yet, this case does not apply to our context.

\biglf
Note this serious difference 
between paths of random fields and functions from
their native spaces. Furthermore, we now consider the random norm of a random function, not random values at certain single points. This is the difference to Section \ref{SecPwReg}.
\biglf
By Theorem \ref{Thegap1}, native space norm convergence in
\eqref{eqpSWR} will fail in general, but there may be weaker norms like $L_2(\Omega)$
that admit convergence of paths in norm. We postpone this to the next section.
\biglf
By Theorem \ref{Thegap1}  we know that IID paths do not lie in
the native space almost surely, but are they almost surely in larger spaces? 
The answer requires a {\em scale} of spaces with different regularity properties. The easiest and more general approach to scaled spaces is via weighted expansions, and we use it for the scale of Sobolev spaces later.
\biglf
To stay close to Fourier series and \MKL{} expansions, we build a scale of Sobolev-type functions over $L_2(\Omega)$ using a fixed basis  of orthonormal functions $v_n$ in $L_2(\Omega)$. 
We penalize the expansion coefficients by positive sequences $\borho=\{\rho_n\}_n$ and define inner products
$$
(f,g)_{L_{2,\borho}(\Omega)}=\sum_{n=1}^\infty (f,v_n)_{L_2(\Omega)}(g,v_n)_{L_2(\Omega)}\rho_n^{-1}.
$$
We mimic the \MKL{} case by assuming $w_n=\sqrt{\lambda_n}v_n$ being orthonormal
in $\calh(C)$ for a covariance function $C$ on $\Omega$,
{see Section \ref{SecOE} for details}.
By comparison of expansions, this implies
$(f,v_n)_{L_2(\Omega)}=\sqrt{\lambda_n}(f,w_n)_{\calh(C)}$, and 
the native Hilbert space for $C$ then consists of all functions $f\in L_2(\Omega)$ with
$$
\sum_{n=1}^\infty (f,w_n)_{\calh(\Omega)}^2=\sum_{n=1}^\infty \lambda_n^{-1}(f,v_n)^2_{L_2(\Omega)}<\infty. 
$$
i.e. it is $L_{2,\bolambda}(\Omega)$.
If $p$ is a path in that space, it is in $L_{2,\borho}(\Omega)$ iff
$$
\sum_{n=1}^\infty (f,v_n)^2_{L_2(\Omega)}\rho_n^{-1}=\sum_{n=1}^\infty \omega_n^2\lambda_n\rho_n^{-1}<\infty.
$$
\begin{Theorem}
For a scale of spaces based on weighted expansions in $L_2(\Omega)$ under the above assumptions,  the native space $L_{2,\bolambda}(\Omega)$ lies in $L_{2,\borho}$ {if} $\sum_{n=1}^\infty \lambda_n\rho_n^{-1}<\infty$. \qed
\end{Theorem}
We now connect this to the scale of Sobolev spaces $W_2^m(\Omega)$ on compact domains $\Omega\subset \R^d$. The variances $\lambda_n^2$
  for \MKL{}
  expansions (\ref{eqMercer}) behave like $n^{-2m/d}$ \citep{santin-schaback:2016-1}
  for $n\to\infty$. 
  To simplify notation, we define an ID path the expansion path \eqref{eqpSWR} for which $S$ is additionally ID. Hence, we have
\begin{Theorem}\label{TheGap2}
If a covariance function generates the native space $W_2^m(\R^d)$, 
and if all other Sobolev spaces are defined via scaling, its ID paths lie almost surely in all Sobolev spaces $W_2^{p}(\R^d)$ for $p<m-d/2$, but almost surely not in $W_2^{m-d/2}(\R^d)$. \qed
\end{Theorem}

On the Sobolev scale, there is a {\em smoothness gap} of order $d/2$ between functions in the native space and paths. This goes back to
\cite{scheuerer:2009-1}, but we provided an explicit and constructive proof
for expansion paths, based on (\ref{eqpSWR}) {and a revised definition
  of paths}, but limited to \MKL{} expansions.
Roughly, if $W_2^m(\R^d)$ is the native space for a covariance kernel,
the borderline space for paths is $W_2^{m-d/2}(\R^d)$ in the above sense.
\biglf
Expansion paths via infinite point sets do not work in the proof,
because we used that the $v_n$ are independent of the smoothness order.
\biglf
It should be  explained why Theorems \ref{TheASPWconv},
\ref{Thegap1} and \ref{TheGap2} are not contradictory.
Imagine a random function generator that produces large numbers of random
paths. When considering only a single point $x$, one gets a random variable
over these paths that has bounded variance.
But something like the $L_2(\Omega)$ norm
of single paths is another random variable over paths that needs a proof
for guaranteeing bounded variance.
\biglf
Theorem \ref{TheGap2} has a positive computational aspect.
If users want to produce cheap random sample {paths} from $W_2^m(\Omega)$ {without any excessive smoothness},
they should not run anything for $n\to \infty$ using {that kernel.
For roughly  the same effect, one can take cases for finite $n$ in
$W_2^{m-d/2}(\Omega)$.}
\biglf
{The paper \citep{steinwart:2019-1} draws much more
detailed conclusions for expansion paths \eqref{eqpSWR}
with respect to various function spaces, but it still starts 
in the standard way from Random Fields. 

{\subsection*{Remarks}
The regularity differences between deterministic interpolation, non-deterministic Kriging, and path construction for random fields should be explained in some more detail, illustrating the above results.}
\begin{enumerate}
    \item {{\bf Deterministic interpolation}}\\
    There is an apparent philosophical discrepancy between the statistics and the numerical analysis approach. For the former, data are realizations from something {\em random}, while for the latter data are the ground truth modulo some additive noise. We start by focusing on the latter: in this case, data are the true values of a function $f$ from the native space of $C$. For finitely many observations from $f$, the smoothness of the (Kriging) interpolant 
is {\em double} with respect to the {smoothness of functions in the native space, because {the interpolant} is generated from  kernel translates. These have the smoothness of the kernel $C$, but native space functions only have the smoothness of the kernel $C_0$ with $C=C_0*C_0$, i.e., $\hat C=(\hat C_0)^2$.}  When the number of observations of $f$ tends to infinity, the {double smoothness is lost by going to the limit in the native space. This well-known {\em smoothness gap} is deterministic and arises when going from translates of $C$ to their Hilbert space limits in $\calh(C)$}. \\
\item {{\bf Kriging}}\\
The variance of BLUP (Kriging) prediction at some target point $x$ is the square of the Power Function, and the latter is the 
kernel conditioned to the data locations. So far, this is only point-dependent, not data-dependent.
The Kriging functions, when extended to all $x$, using data of some $f$ at the nodes, {coincide with deterministic interpolation}. Yet, the assumption that {\em ground truth} means that data come from a function
in the native space is a questionable  addition. The "escape scenario" {\citep{karvonen:2023-1}} results of Numerical Analysis 
show that if data come from a function in a larger Sobolev space, the Kriging solution is still
convergent to the data function, albeit in the weaker Sobolev norm {\citep{narcowich-et-al:2005-1,narcowich-et-al:2006-1}}. This is often called 
``misspecification'' \citep{stein:1999-1}, because the $C$ model is different from the $f$ model. 
For a finite number of observations, the Kriging solution has the smoothness of $C$-translates. 
It is roughly twice the $\calh(C)$ smoothness, {as pointed out above.
\item {\bf{Paths of Random Fields}}} \\
{Here, there is no given function $f$ and no given values at points}.
Using, {\em e.g.}, {the Cholesky decomposition of the kernel matrix and a partial sum of (\ref{eqpSWR})}, one can generate $n$ random data values 
at each point $x$ of the domain. {By (\ref{eqresvar})},   the variance for random selection of a new value will get smaller when increasing $n$.
After $n$ steps, one has  $n$ random data values that will be the 
final path values. If these values are taken as ``ground truth'' {like Interpolation or Kriging}, they have a Kriging interpolant that lies in 
the native space and even has the excess smoothness of kernel translates. But there is no $f$ that supplied them,
and the ``noise-generating'' process of adding new random sample values, though using smaller variances
from step to step, induces more smoothness loss than {the deterministic loss described above.}
It is exactly this additional smoothness loss by randomness that forms the ``gap'' {described by Theorem \ref{TheGap2}. 
\item Summarizing, it can be seen as a miracle} that the randomness of the new samples 
is limited by $C$ in such a way {(namely by 
(\ref{eqresvar})}) that one safely arrives in larger Sobolev spaces {instead of at what Numerical Analysts would call noise, i.e. a function with hardly any regularity}.
Yet, the {\em exact} Sobolev limit space cannot be reached, {by Theorem \ref{TheGap2}}. 
\end{enumerate}

  \section{Discussion and Open Problems}
This section serves both as a rejoinder as well as an introduction to open problems. We sketch them below. 
\begin{itemize}
\item Turning the standard approach to paths upside down 
  avoids all problems with joint probability spaces on infinite sets,
  and it is closer to {\em Nature} because it goes from global to local. 
\item {We consider the approach proposed in this paper} more transparent {with respect to earlier literature}:  the random variables
  $R_C(x)$ defined via $C$ are just the values of the paths at points $x$.
  This is possible because paths are defined before random variables
  are defined. 
\item It simplifies the analysis of paths by the explicit form
  \eqref{eqpSWR} and the detour via the Hilbert space $\calh(C)$.
\item Since the bottom-up path construction in the Gaussian case
  is a special case, nothing is lost in that case.
\item {There is a substantial novelty in that}  now other choices for introducing randomness after picking a
  covariance kernel are possible. 
\item 
The standard sampling algorithms {based on infinite dense point sets} are fast if $C$ is smooth, {because the residual variance 
$C_{N+1}(x,x)$ from \eqref{eqCCww} 
decays quickly, leading to a fast convergence in $L_2(\Omega)$ of partial sums \eqref{eqpSWR} of paths. These are in the native space,
  i.e. somewhat } too smooth, unless a kernel for a larger native space is used.\\
  {Therefore the literature has various techniques with better numerical
  complexity for low smoothness, A widely used method based on \citep{lindgren-et-al:2011-1} uses weak solutions of stochastic differential equations. Another strategy, starting from \cite{vecchia:1988-1} exploits  that inverses of kernel matrices
  have good approximations by sparse matrices. These techniques calculate approximations to paths, but not paths themselves, and it should be checked if our approach to paths makes the error analysis easier. } 
\item It needs further work to study the statistical differences
  of path calculations with 
  different
  univariate probability distributions $\calr$. Similarly,
  the aforementioned methods for calculating path approximations
  should be compared.
\item Since this paper proceeds via the Hilbert space $\calh(C)$,
  an extension to the vector-valued case should be possible using
  the approach to paths used here.
\end{itemize}
\bibliographystyle{apalike}
\bibliography{bib_matricial, RSbib}
\end{document}